\documentstyle{article}
 \def\vt{t\kern-0.22em\raise.18ex\hbox{\char'47}\lower.18ex\hbox{}\kern-0.08em}
 
\newtheorem{th}{Theorem}[section]

\newtheorem{lm}{Lemma}[section]

\newtheorem{con}{Conjecture}[section]

\newtheorem{rem}{Remark}[section]

\newcommand{\old}[1]{{}} 
\newcounter{obr}
\newcounter{tabul}
\begin{document}
\title{ On the largest reduced neighborhood clique cover number  of a graph 
}
\author{Farhad Shahrokhi\\
Department of Computer Science and Engineering,  UNT\\
farhad@cs.unt.edu
}

\date{}
\maketitle
\thispagestyle{empty}
\date{} \maketitle
\begin{abstract}  Let  $G$ be a graph and $t\ge 0$.  
A new graph parameter termed the largest reduced neighborhood clique cover number of $G$, denoted by 
${\hat\beta}_t(G)$, 
is introduced.  Specifically,  ${\hat\beta}_t(G)$ is the largest, overall $t$-shallow minors $H$ of $G$, of the smallest number of  cliques that can  cover   any  closed neighborhood of a vertex in $H$.  We verify that  ${\hat\beta}_t(G)=1$ when $G$ is chordal, and,  ${\hat\beta}_t(G)\le  s$, where $G$ is  an incomparability graph that does not have a $t-$shallow minor which is isomorphic to an induced star on $s$ leaves. 
Moreover, general  properties of  ${\hat\beta}_t(G)$ including the connections to  the greatest reduced average 
density of $G$, or $\bigtriangledown_t(G)$ are studied and investigated. For instance we show ${{\hat\beta}_t(G)\over 2}\le \bigtriangledown_t(G)\le p.{\hat\beta}_t(G),$ where 
$p$ is the size of a largest complete graph which is a $t-minor$ of $G$. Additionally we prove that 
   largest ratio of any minimum clique cover to the maximum independent set taken overall $t-$minors of $G$ is a lower bound for ${\hat\beta}_t(G)$. 
We further introduce the class of bounded neighborhood clique cover number  for which ${\hat\beta}_t(G)$ has a finite value for each $t\ge 0$ and verify the membership of  geometric intersection graphs of fat objects (with no restrictions on the depth)  to this class. The results support the  conjecture that the class graphs with polynomial bounded neighborhood clique cover number  may have separator theorems with respect to certain measures. 
 \end{abstract}

\section{Introduction}
  
Throughout this paper $G$ denotes an undirected graph on vertex set $V(G)$ and the edge set $E(G)$.
Let $\delta(G)$ denote the minimum degree of $G$. Recall that  the degeneracy of $G$, denoted by ${\hat\delta}(G)$ is the largest minimum degree among all induced subgraphs of $G$.   A  {\it $t-$shallow minor}, or a {\it $t-$minor} of $G$ in short,  is a minor   of $G$ which is obtained by contracting  connected subgraphs of radius  at most $t$, and deleting vertices (but not edges).

Shallow minors were  first introduced by Leiserson and was used as a tool by Plotkin et al. \cite{P} to investigate the 
graph separators and gap between certain Multi-flows and ratio cuts. 

 Nesetril and Ossona de Mendez  introduced the fundamental  notion of $\bigtriangledown_t(G)$, termed the greatest reduced average density of $G$ (grad of $G$ in short),  as well as, the class of  bounded expansion graphs and have shown many  important properties and applications for this class \cite{N1,N2,N3,N4}. More precisely, they defined  $\bigtriangledown_t(G)$ to be the maximum edge density of any $t-$minor in $G$. 
  
 
 Nesetril and Ossona de Mendez and define $G$ to have bounded expansion, 
 if
$\bigtriangledown_t(G)$  
is finite for every $t\ge 0$.  This class of graphs contains many traditionally known ``sparse'' graphs. Additionally,
there is an intimate connection between bounded expansion and graph separators which in return gives rise to nice algorithmic and structural results.

 The motivation in writing this paper arises from the fact that although $\bigtriangledown_t(G)$
   is very  helpful in studying  the properties of ``sparse'' graphs,  what if  the  underlying graph is not``sparse''?

 For a graph $H$,  let $\beta(H)$ denote the clique cover number of $H$, that is, the minimum number of cliques that partition $V(H)$.   
 Now  for any $x\in V(H)$,  let $H_x$ denote the 
 the closed neighborhood of $x$ in $H$, and let ${\tilde \beta(H)}=min_{x\in V(H)}\{\beta(H)\}$.
   For any graph $G$ and $t\ge 0$  define ${\hat\beta}_t(G)$ to be  the largest value of $\tilde\beta(H)$ for any $t-$minor $H$ of $G$.  We call ${\hat\beta}_t(G)$ the largest  reduced neighborhood clique cover number  of $G$. We say $G$ has bounded  neighborhood clique cover number  
   if ${\hat\beta}_t(G)$ has finite value for each $t\ge 0$. Note that $\hat\beta_t({K_n})={1}$  for any $t\ge 0$, nonetheless 
 $\bigtriangledown_t(K_n)={n-1\over 2}$.   
 
   In Section Two we investigate the largest reduced neighborhood clique cover number  in chordal and incomparability graphs. 
   We show$\hat\beta_t(G)=1$ for any chordal graph $G$, whereas $\hat\beta_t(G)\le s_t$, for any incomparability graph $G$,  where $s_t$ is the size (number of leaves) in a  largest $t-minor$ which is isomorphic to an induced  star.
  Furthermore we  show  ${{\hat\beta}_t(G)\over 2}\le \bigtriangledown_t(G)\le p.{\hat\beta}_t(G),$ where 
  $p$ is the size of a largest complete graph which is a $t-minor$ of $G$. Additionally we prove that 
   largest ratio of any minimum clique cover to the maximum independent set taken overall $t-$minors of $G$ is a lower bound for ${\hat\beta}_t(G)$. In Section Three we verify that the intersection graph of fat objects has a bounded neighborhood
  clique cover number (in fact polynomial),  provided that the geometric dimension is bounded. In light of existence  separation theorems for fat objects in bounded dimension, we conjecture that 
  graphs with polynomial bounded neighborhood clique cover number have separation theorem with respect to a measure.
  (See Conjecture \ref{co1}.)  
  
   \section{Main Results}
   Recall that a chordal graph does not have any chord-less cycles \cite{Go}. An incomparability graph is a graph whose complement has a transitive orientation \cite{Tr}.
   
   \begin{th}\label{t1}
   
   {\sl ~~~~~~~~~~~~~~~~~~~~~~~~~~~~~~~~~~~~~~~~~~~~~~~~~~~\\
   
   $(i)~$Let $G$ be a chordal graph, then ${\hat\beta}_t(G)=1$. 
   
   $(ii)~$Let $G$ be an incomparability graph that does not have a $t-$shallow minor which is isomorphic to an induced star on $s$ leaves.  Then, ${\hat\beta}_t(G)\le  s$.

   }
   \end{th}
   
   {\bf Proof}.   For $(i)$ observe that any chordal graph has a simplicial vertex, that is a vertex whose neighborhood is clique \cite{Go}. Now observe that  the  class of chordal graphs is closed under contraction and deletion of vertices. So any minor of $G$, say $H$, which is obtained from $G$ by contraction and deletion of vertices is also chordal.  Since $H$ is chordal it has vertex $x$, so that $H_x$ is a clique  \cite{Go}.  Thus, $\beta(H_{x})=1$, implying ${\tilde\beta}(H)=1$, and consequently ${\hat\beta}_t(G)=1$ as claimed. 
   
For $(ii)~$, let  $H$ a $t-$minor of $G$, and note that $H$ is an incomparability graph. Now let $C=\{C_0, C_1,...,C_k\}$ be a greedy
clique cover in $H$, where $k+1$ is the size of  a largest  independent set in $H$.
Thus, $C_1$ is the set of all sources
in the transitive orientation of complement of $H$, 
$C_2$ is the set of all sources that are obtained  after removal of $C_1$,
etc. 
Now let $a\in C_0$.  Let $N(i)$ denote the set of all vertices adjacent to $a$ which are in $C_i$ and   
let $q$ be the largest $i$ so that there is a vertex in $C_i$ adjacent to $x_0, i=0,1,...,k$.   

Now let  $x_q\in C_q$,  then for $i=q-1,...,0$ there
is $x_i\in C_i$ so that  $x_ix_{i+1}\notin E( H)$. It follows that
for $i=0,1,...,q$, we have $ax_i\in E({H})$ and thus for $j>i$, we have 
$x_ix_j\notin E(G)$. We conclude $a$ has at least one neighbor in each $N(i),i=0,1,2,...,q$, and that the induced graph on $x_0,x_1,x_{i+1},...,x_q$
is a star in $H$, with center $a$, having $q$ leaves. Consequently $\beta(H_{x_0})= q+1\le s+1$ and hence the claim follows.

   $\Box$
   
  It is convenient to   define $\hat\bigtriangledown_t(G)$ denote the largest  minimum degree  of any $t-$minor in $G$.
Observe that $\bigtriangledown_t(G)\le \hat\bigtriangledown_t(G)\le 2\bigtriangledown_t(G)$.

 \begin{th}\label{t2}
 {\sl Let $t\ge 0$, then for any graph $G$  
 
$${{\hat\beta}_t(G)\over 2}\le \bigtriangledown_t(G)\le p.{\hat\beta}_t(G),$$

where $p$ is the size of a largest $t-$minor in $G$ which is isomorphic to $K_p$. 
}
\end{th}

{\bf Proof. } For the lower bound, let $H$ be a $t-$minor of $G$, then for any $x\in V(G)$, we must have 
$\beta(H_x)\le deg_H(x)$,  where $deg_H(x)$ is degree of $x$ in $H$,  and hence $\tilde\beta(H)\le \delta(H)$. 
Consequently,  $\hat\beta(G)\le \hat\bigtriangledown_t(G)\le 2\bigtriangledown_t(G)$. 
For the upper bound, let   $H$ be  a $t-$minor of $G$, then, for any $x\in V(G)$, we have  $deg_H(x)\le p.\beta(H_x)$. 
Consequently, for any $H$ which is a $t-$minor of $G$, we have  $\delta(H)\le p.{\tilde\beta (H)}$.  
Hence, $\hat\bigtriangledown_t(G\le p.\hat\beta_t(G) $ and claim follows since   $\bigtriangledown_t(G)\le \hat\bigtriangledown_t(G)$.
$\Box$\

\begin{th}\label{t3}
{\sl Let $G$ be a graph, let $t\ge 0$ and let $H$ be a $t-$minor of $G$. Then 
 ${\beta(H)\over \alpha(H)}\le \hat\beta_t(G)$, where $\alpha(H)$ is the size of a largest independent set in $H$.
}
\end{th}

{\bf Proof.} We apply a modified  version of the algorithm due to  Erfat et al \cite{EKNS} for approximating the maximum independent set  of fat objects.  Let$H^0=H$, and pick a vertex  $x_0\in V(H^0)$ with 
$\beta(H^0_{x_0})=\tilde\beta(H^0) $. Now for $i\ge 1$ (at iteration $i$ of algorithm), let $H^i$ be the subgraph of $H$ on the vertex set $V(H^i)=
V(H)-\cup_{j=0}^{i-1}V(H^j_{x_j})$, note that $H^i$ is a$t-$minor of $G$, and pick  $x_i\in V(H^i)$ with $\beta(H^i_{x_i})=\tilde\beta(H^i)$. Upon termination of the algorithm, we obtain a set of independent vertices in $H$, 
$S=\{x_0,x_1,....,x_r\}$, and a set of $B=\{B_0,B_1,...,B_r\}$ where for $i=0,1,..., r$, $B_i$ is a clique cover
for $H^i_{x_i}$ with  $|B_i|=\tilde\beta(H^i)\le {\hat\beta}_t(G)$. Thus we can cover $H$ with at most 
$(r+1).{\hat\beta}_t(G)$ cliques which implies the claim.  $\Box$.
\section{Connections to Geometry }
The  size of an object in ${R}^d$ is the side length of its smallest 
enclosing hypercube. 
Several definitions for fat  objects exists. We use the following from \cite{Chan}.
Let ${\cal C}$ be a collection  of objects in ${R}^d$, and let
$c\ge 1$.
The objects in ${\cal C}$ are $c-${\it fat},   
if  for any $r\ge 0$, and any   Box  $B$ of size $r$, there is a set of points
in ${R}^d$ of cardinality at most $c$,  so that 
any object of size at least $r$ in ${\cal C}$ that intersects $B$ 
contains one of these points.  Note  that  spheres,  cubes, and boxes with bounded aspect ratios, or  ``suitable mixture'' of these objects are $c-$fat objects, for different values of $c$. Typical values of $c$ is usually $O(b^d)$, for some constant $b\ge 2$.

\begin{th}\label{t6}
{\sl Let ${\cal C}$ a  collection  of $c-$fat objects in ${R}^d$,  let $t\ge 0$,  and let $G=G({\cal C})$ denote the intersection graph of ${\cal C}$. 
Then ${\hat \beta}_t(G)=O(c.t^{2d})$. 
}
\end{th}
To prove the claim we need the following Lemma which may be folklore in the community. This author verified it  in a unpublished manuscript in 2013. For the  published version of a similar result  see the recent work of Har-Peled and Quanrud \cite{HQ}.
\begin{lm}
{\sl
Let ${\cal C}$ be a  collection  of $c-$fat objects in ${R}^d$, and let 
${\cal C}^t$ be a collection  of 
subsets of ${\cal C}$,  with the property that, for each subset  ${\cal X}$ in ${\cal C}^t$,  the  diameter of the induced subgraph of $G$ of ${\cal X}$, denoted by  $G_{\cal X}$,  is bounded by $t$. Then, ${\cal  C}^t$ is a
set of $c'-$fat objects in ${R}^d$, where, $c'=O(c.t^{2d})$. 
}

\end{lm}

{\bf Proof Sketch.}  By the above Lemma any $t-$minor $H$ of $G$ is the intersection 
graph of $c'-$fat objects with $c'=O(c.t^{2d})$. It further follows from the definitions given above that
for any $t-$minor $H$,  ${\tilde\beta}(H)=O(c.t^{2d})$ and hence ${\hat\beta}_t(H)=O(c.t^{2d})$. $\Box$.

\begin{rem}
{\sl It was observed  in \cite{HQ} that the intersection  graph $G$  of a set of  $c-$fat objects has bounded 
expansion, provided that the depth (the maximum number objects that contains 
any point), $d$, and $c$ are bounded. This also follows from Theorems 2.2 and 3.1. 
Combining this with a result of Dvorak \cite{D2}, one can show (under the above  assumptions)  that  for any integer $k$,
the minimum number of vertices in $G$ that dominate all vertices at   distance at most $k$,
is bounded above by  a linear function of maximum number of vertices whose  pairwise distances are at least $2k+1$. (See \cite{HHP} for general domination theory, \cite{HOS} for distance  domination, and  specifically  \cite{BM} for a related result on distance domination  when $G$ excludes a large bipartite minor.}) 
\end{rem}

Alber and Fiala \cite{AF} studied independent and dominating sets in the intersection graph of two dimensional discs, and derived separation theorems for these graphs.  Chan \cite{Chan} derived generalized separation theorems for fat objects with respect to certain ``measures''. (Note that since intersection graph of fat objects may be very dense,  meaningful separation theorems with respect to number of vertices do not exist.)  Motivated by the earlier results of Chan \cite{Chan} and Alber and Fiala 
\cite{AF} on geometric separation of 
fat objects (with respect to a measure),  we have  derived in \cite{Sh1} 
a general combinatorial separation theorem 
whose statement without a proof
was announced in \cite{Sh2}.  In addition,  we provided a refinement of the separation theorems  of Chan in \cite{Sh3} which  led to the first sub exponential time  algorithms for computing the independent sets in fact objects. In light of existence of separation theorems for fat objects, and since  for fat objects the parameter $c$, is bounded by a function of dimension $d$, and since  by Theorem \ref{t6}, the parameter ${\hat\beta}_t(G)$ ($G$ being intersection graph) is polynomially bounded by $t$, we suspect that suitable separation theorems may exist for classes of graphs that have a 
(polynomially)  bounded  neighborhood clique cover number.  

\begin{con}\label{co1}
{\sl Let  $t\ge 0$, then the class of all graphs $G$, for which $\hat\beta_t(G)$ is polynomially bounded by  $t$
has a ``suitable'' separation theorem with respect to the measure  $\alpha$. That is, there is $S\subseteq V$
with $\alpha(S)=O({\alpha(G)}^{1-\epsilon})$, so that for each component  $C$ of $G-S$, we have 
$\alpha(C)\le {\alpha(G)\over 2}$. Here  for any $A\subseteq V$, $\alpha(A)$ denotes the size of a largest independent set in the induced subgraph on $A$.

}

\end{con}

We remark that the above conjecture is supported by the evidence that polynomially bounded expansion graphs
have sub-linear separators with respect to number of vertices \cite{D1}.

\end{document}